\documentclass[a4paper,11pt]{article}
\usepackage{amscd, amsmath, amssymb, amsfonts, epic}
\usepackage{graphicx}
\usepackage{a4wide}
\usepackage{xypic}
\usepackage{paralist}
\usepackage{booktabs}
\usepackage{multirow}

\title{Finiteness results for 3-folds with semiample anticanonical
  bundle} \author{Artie Prendergast-Smith}

\date{}

\def\Z{\text{\bf Z}}
\def\Q{\text{\bf Q}}
\def\R{\text{\bf R}}

\def\P{\text{\bf P}}

\def\arrow{\rightarrow}

\def\iso{\cong}
\def\curlyo{\mathcal{O}}
\def\Pic{\text{Pic}} 
\def\Aut{\text{Aut}}
\def\PsAut{\text{PsAut}}
\def\Div{\text{Div}}

\newcommand{\Nef}[1]{\overline{A(#1)}}
\newcommand{\Mov}[1]{\overline{M(#1)}}
\newcommand{\Curv}[1]{\overline{\text{Curv}(#1)}}

\newtheorem{theorem}{Theorem}[section]
\newtheorem{lemma}[theorem]{Lemma}
\newtheorem{corollary}[theorem]{Corollary} 
\newtheorem{proposition}[theorem]{Proposition}

\newtheorem{conjecture}[theorem]{Conjecture} 

\newtheorem{examples}[theorem]{Example}

\begin{document}
\maketitle

The purpose of this paper is to give some evidence for the
Morrison--Kawamata cone conjecture for klt pairs. Roughly speaking,
the cone conjecture predicts that in appropriate `Calabi--Yau-type'
situations, the groups of automorphisms and pseudo-automorphisms of a
projective variety act with rational polyhedral fundamental domain on
the nef and movable cones of the variety. (See Section
\ref{section-conj} for the precise statement.)

In this paper we prove some statements in this direction, in the case
of a mildly singular 3-fold with semiample anticanonical bundle of
positive Iitaka dimension. Let us say a bit about how such a 3-fold
looks geometrically. The anticanonical bundle defines a contraction
morphism $X \arrow S$ to a positive-dimensional base; by adjunction
all smooth fibres are varieties whose canonical bundle is torsion. So
the generic fibre is a point, an elliptic curve, or a Calabi--Yau
surface, according as the Iitaka dimension is 3, 2, or 1. (Here a
Calabi--Yau surface means an abelian, K3, Enriques or hyperelliptic
surface.) If the contraction morphism is equidimensional,
classification results due to Kodaira and Miranda (for fibre dimension
1) and Kulikov and Crauder--Morrison (for fibre dimension 2) give
information about the singular fibres. (See
\cite{FriedmanMorrison1981} for details of these classification
results.)

We will see that all 3-folds of this kind fall inside the scope of the
Morrison--Kawamata cone conjecture. Our main result is the following
finiteness theorem, which can be regarded as a weak form of the
conjecture for these varieties.
\begin{theorem} \label{theorem-maintheorem1}
Let $X$ be a $\Q$-factorial terminal Gorenstein $3$-fold with $-K_X$
semiample of positive Iitaka dimension. Let $f: X \arrow S$ be the
contraction morphism given by sections of some power of $-K_X$. Then
the effective movable cone $\Mov{X}^e$ decomposes as the union of the
effective nef cones of the small $\Q$-factorial modifications of $X$,
all of these small $\Q$-factorial modifications are over $S$, and the
decomposition is finite up to the action of the group $\PsAut(X/S)$ of
pseudo-automorphisms of $X$ over $S$.
\end{theorem} 
For terminal 3-folds, Gorenstein singularities are the same thing as
hypersurface singularities. The latter description is geometrically
clearer, but we use the former description to remain consistent with
the literature. In any case, we make the Gorenstein assumption so that
we can use a theorem of Mori and its extension by Cutkosky classifying
certain contraction morphisms on terminal 3-folds.

We will explain in Section \ref{section-conj} why Theorem
\ref{theorem-maintheorem1} can be regarded as a weak form of the cone
conjecture for the given class of varieties. The key ingredients in
the proof of the theorem are Kawamata's corresponding result for the
relative movable cone $\Mov{X/S}^e$ and the Mori--Cutkosky
classification result mentioned above. Indeed, given these two
results, our theorem follows as a simple consequence. The proof is
given in Section \ref{section-mainproof}.

In Section \ref{section-movablecone} we study the movable cone in more
detail, under the stronger assumption that $-K_X$ has Iitaka dimension
2 and the contraction $X \arrow S$ is equidimensional. In particular
we get a complete description of the relative movable cone (Theorem
\ref{theorem-relativemovable}), which leads to a finiteness result for
the (absolute) movable cone (Theorem \ref{theorem-cone}) and a proof
of the full movable cone conjecture in some examples (Example
\ref{example-quadrics}).

In Section \ref{section-nefcone} we give some partial results
concerning the conjecture on the nef cone. In particular, for the
simple case of smooth non-rationally connected $X$ such that $-K_X$ is
semiample of Iitaka dimension 2, we show how the nef cone conjecture
follows from classification results of Bauer--Peternell (Theorem
\ref{corollary-nefcone}).

Thanks to Burt Totaro for many helpful comments.

\section{The cone conjecture} \label{section-conj}

We work throughout over an algebraically closed field $k$ of
characteristic 0.

Suppose $f: X \arrow Z$ is a projective surjective morphism of normal
varieties with connected fibres. A Cartier divisor $D$ on $X$ is said
to be {\it nef over $Z$} (resp. {\it movable over $Z$}, {\it big over
  $Z$}, {\it effective over $Z$}) if $D \cdot C \geq 0$ for all curves
$C$ mapped to a point by $f$ (resp.  if $\operatorname{codim \ Supp
  \ Coker } (f^*f_* \mathcal{O}_X(D) \arrow \mathcal{O}_X(D)) \geq 2$,
if the Iitaka dimension $\kappa(X_\eta, D_\eta) = \operatorname{dim }X
- \operatorname{dim } Z$ for $\eta \in Z$ the generic point, if $f_*
\mathcal{O}_X(D) \neq 0$).

We define the real vector space $N^1(X/Z)$ to be
$\left(\Div(X)/\equiv_Z \right) \otimes \ \R$ where $\Div(X)$ is the
group of Cartier divisors on $X$ and $\equiv_Z$ denotes numerical
equivalence over $Z$. The {\it relative nef} cone $\Nef{X/Z}$
(resp. {\it closed relative movable cone} $\Mov{X/Z}$) is the closed
convex cone generated by classes of Cartier divisors which are nef
over $Z$ (resp. movable over $Z$). The {\it relative big} cone
$B(X/Z)$ is the open cone generated by Cartier divisors which are big
over $Z$. The {\it relative effective cone} $B^e(X/Z)$ is the cone
generated by Cartier divisors which are effective over $Z$. We denote
by $\Nef{X/Z}^e$ and $\Mov{X/Z}^e$ the intersections $\Nef{X/Z} \cap
B^e(X/Z)$ and $M(X/Z) \cap B^e(X/Z)$, and call them the {\it relative
  effective nef cone} and {\it relative effective movable cone}
respectively.

For later use, we also make the following definitions. The real vector
space $N_1(X/Z)$ is defined to be the dual of $N^1(X/Z)$;
equivalently, it is the real vector space spanned by numerical classes
of curves on $X$ which map to a point on $Z$. We define the {\it
  relative closed cone of curves}, denoted $\Curv{X/Z}$, to be the
closed convex cone in $N_1(X/Z)$ dual to the cone
$\Nef{X/Z}$. (Equivalently, it is the closed convex cone spanned by
the classes of irreducible curves on $X$ which map to a point on $Z$.)
We should stress that our notation for the cone of curves is
nonstandard: it replaces the standard notation
$\overline{NE(X/Z)}$. The new notation is intended to be more
informative.

For all the notation introduced above, in the case that $Z=
\operatorname{Spec }k$ we will omit it from the notation and simply
write $N^1(X)$, $\Nef{X}$, and so on. In this case we also omit the
adjective `relative' from all the corresponding terms. 

Define a {\it pseudo-isomorphism} from $X_1$ to $X_2$ over $Z$ to be a
birational map $X_1 \dashrightarrow X_2$ over $Z$ which is an
isomorphism in codimension 1. For a $\Q$-factorial variety $X$ over
$Z$, a {\it small $\Q$-factorial modification} (SQM) of $X$ over $Z$
means a pseudo-isomorphism over $Z$ from $X$ to another $\Q$-factorial
variety with a projective morphism to $Z$. Note that if $\alpha: X'
\dashrightarrow X$ is an SQM over $Z$, there is a canonical
identification $N^1(X'/Z) \iso N^1(X/Z)$ given by proper transform of
divisors. This identification maps the effective movable cone of $X'$
to that of $X$, so in particular identifies the effective nef cone
$\Nef{X'/Z}^e$ with a subcone of $\Mov{X/Z}^e$. This will explain the
second statement of Conjecture \ref{conj-totaro}. It is important to
note that the identification depends on the map $\alpha$ and not just
the variety $X'$; we will write $\Nef{X'/Z, \alpha}$ for the image of
the nef cone inside $N^1(X/Z)$ when it is important to keep track of
the identification.

For an $\R$-divisor $\Delta$ on a $\Q$-factorial variety $X$, the pair
$(X,\Delta)$ is {\it klt} if for all resolutions $\pi: \tilde{X}
\arrow X$ with a simple normal crossing $\R$-divisor $\tilde{\Delta}$
such that $K_{\tilde{X}}+\tilde{\Delta} = \pi^*(K_X+\Delta)$, the
coefficients of $\tilde{\Delta}$ are less than 1. (For instance if $X$
is smooth and $D$ is a smooth divisor on $X$, then $(X,rD)$ is klt for
any $r<1$.) We say that $(X/Z,\Delta)$ is a {\it klt Calabi--Yau pair}
if $(X,\Delta)$ is a $\Q$-factorial klt pair with $\Delta$ effective
such that $K_X+\Delta$ is numerically trivial over $Z$.

We denote the groups of automorphisms or pseudo-automorphisms of $X$
over $Z$ which preserve a divisor $\Delta$ by $\Aut(X/Z,\Delta)$ and
$\PsAut(X/Z,\Delta)$. A {\it rational polyhedral cone} in $N^1(X/Z)$
means a convex cone spanned by a finite set of classes of Cartier
divisors.

\begin{conjecture} \label{conj-totaro}
Let $(X/Z, \Delta)$ be a klt Calabi--Yau pair. Then:

\quad (1) The number of $\Aut(X/Z,\Delta)$-equivalence classes of
faces of the relative effective nef cone $\Nef{X/Z}^e$ corresponding
to birational contractions or fibre space structures is
finite. Moreover, there exists a rational polyhedral cone $\Pi$ which
is a fundamental domain for the action of $\Aut(X/Z,\Delta)$ on
$\Nef{X/Z}^e$ in the sense that

\quad (a) $\Nef{X/Z}^e = \Aut(X/Z,\Delta) \cdot \Pi$,

\quad (b) $\text{Int} \ \Pi \cap g_* \text{Int} \ \Pi = \emptyset$ for
any $g \in \Aut(X/Z,\Delta)$ such that $g_* \neq 1$ (where $g_* \in
GL(N^1(X/Z))$ is the linear automorphism induced by $g$).

\quad (2) The number of $\PsAut(X/Z,\Delta)$-equivalence classes of
chambers $\Nef{X'/Z,\alpha}^e$ in the relative effective movable cone
$\Mov{X/Z}^e$ corresponding to SQMs $\alpha: X' \dashrightarrow X$ of
$X$ over $Z$ is finite. Moreover, there exists a rational polyhedral
cone $\Pi'$ which is a fundamental domain for the action of
$\PsAut(X/Z,\Delta)$ on $\Mov{X/Z}^e$.
\end{conjecture}
The first statement in parts (1) and (2) of the conjecture follows
from the second.  We refer to the first statement in parts (1) and (2)
as the {\it weak (nef or movable) cone conjecture}.

Now we explain how the conjecture relates to Theorem
\ref{theorem-maintheorem1}. Suppose $X$ is a $\Q$-factorial terminal
Gorenstein 3-fold with $-K_X$ {\it semiample} (meaning that for some
positive integer $m$ the line bundle $-mK_X$ is basepoint-free) and of
positive Iitaka dimension. Since $X$ is terminal Gorenstein, it has
isolated singularities \cite[Theorem 1.1]{Reid1983}. Therefore
choosing $m$ large enough so that $-mK_X$ is basepoint-free, a general
divisor $D \in |-mK_X|$ does not intersect the singular locus of
$X$. Moreover, Bertini's theorem says that a general $D \in |-mK_X|$
is smooth outside the singular locus of $X$, so by the previous
sentence such a divisor $D$ must be smooth. Any resolution of $X$ is
then a log resolution of $(X,D)$, and on such a resolution the
discrepancies of all exceptional divisors are the same for $(X,D)$ as
for $(X,0)$. Since $X$ is terminal, these are in particular
positive. Assuming without loss of generality that $m>1$, and putting
$\Delta = \frac{1}{m}D$, the pair $(X,\Delta)$ is therefore a klt pair
with $\Delta$ effective, and $K_X + \Delta = 0$ in $N^1(X)$. In other
words $(X,\Delta)$ is a klt Calabi--Yau pair (over $Z=
\operatorname{Spec } k$). So the weak movable cone conjecture predicts
that there should be finitely many $\PsAut(X,\Delta)$-orbits of nef
cones inside the effective movable cone $\Mov{X}^e$. Finally since we
have a morphism $f:X \arrow S$ given by a multiple of the line bundle
$-K_X$, the divisor $\Delta$ is a pullback of a $\Q$-divisor on
$S$. Therefore the group $\PsAut(X/S)$ is a subgroup of
$\PsAut(X,\Delta)$. Theorem \ref{theorem-maintheorem1} then implies
that the weak movable cone conjecture holds for $(X,\Delta)$.

We conclude this section by outlining the history and current status
of Conjecture \ref{conj-totaro}. Inspired by mirror symmetry, Morrison
\cite{Morrison1992} first proposed the conjecture for Calabi--Yau
varieties. This was generalised to Calabi--Yau fibre spaces in
\cite{Kawamata1997} and to klt Calabi--Yau pairs in
\cite{Totaro2010}. The conjecture was proved for Calabi--Yau surfaces
by Sterk--Looijenga, Namikawa, and Kawamata
\cite{Sterk1985,Namikawa1985,Kawamata1997}, for 3-dimensional
Calabi--Yau fibre spaces over a positive-dimensional base by Kawamata
\cite{Kawamata1997}, and for klt Calabi--Yau pairs of dimension 2 by
Totaro \cite{Totaro2010}. In dimension 3 the conjecture remains open,
although there are significant results due to Oguiso--Peternell
\cite{Oguiso2001}, Szendr\"oi \cite{Szendroi1999}, Uehara
\cite{Uehara2004}, and Wilson \cite{Wilson1992}. Finally, there are
verifications of special cases in dimension 3 such as fibre products
of rational elliptic surfaces by Grassi--Morrison
\cite{GrassiMorrison1993} and Horrocks--Mumford quintics by Borcea
\cite{Borcea1991} and Fryers \cite{Fryers1999}.

\section{Proof of Theorem \ref{theorem-maintheorem1}} \label{section-mainproof}

In this section we give a proof of Theorem \ref{theorem-maintheorem1},
thereby showing that the weak movable cone conjecture holds for
$\Q$-factorial terminal Gorenstein 3-folds with anticanonical
bundle semiample of positive Iitaka dimension. As mentioned in the
introduction, the proof follows easily from Kawamata's corresponding
result in the relative case and Mori's classification of $K$-negative
extremal 3-fold contractions. The basic idea is to prove that any
non-nef effective movable divisor on our 3-fold can be made nef by
a sequence of flops, and that these flops are compatible with the
morphism $f: X \arrow S$.

The first step is to use the the classification of $K$-negative
extremal 3-fold contractions, due to Mori \cite{Mori1982} and
extended to the singular case by Cutkosky \cite{Cutkosky1988}. Here a
{\it Mori fibre space} means a contraction morphism $f: X \arrow Z$
with $\operatorname{dim } Z< \operatorname{dim }X$, relative Picard
number $\rho(X/Z) =1$, and $-K_X$ ample over $Z$. (For later use, we
give the classification statement in more detail than we currently
need.)

\begin{theorem}[Mori--Cutkosky] \label{theorem-moriclassification}
Suppose $X$ is a $\Q$-factorial terminal Gorenstein 3-fold. Suppose
$R \subset \Curv{X}$ is $K$-negative extremal ray and $f: X \arrow Z$
is the contraction of $R$. Then either $\operatorname{dim } Z \leq 2$
and $f: X \arrow Z$ is a Mori fibre space, or else $f$ is birational,
the exceptional set $\text{Exc}(f)$ is a prime divisor $D$ on $X$, and
the possibilities for $D$ and $f$ are as follows:
\begin{enumerate}
\item $f:X \arrow Z$ is the blowup of the ideal sheaf $I_C$ of a
  reduced irreducible lci curve $C \subset Z$ with exceptional divisor
  $D$; 
\item $D \iso \P^2$ with normal bundle $\mathcal{O}_D (D) \iso
  \mathcal{O}_{\P^2}(-1)$, and $f$ contracts $D$ to a smooth point;
\item $D \iso \P^1 \times \P^1$ with $\mathcal{O}_D (D)$ of bidegree
  $(-1,-1)$, and $f$ contracts $D$ to a point;
\item $D$ is isomorphic to a singular quadric in $\P^3$ with
  $\mathcal{O}_D (D) = \mathcal{O}_D \otimes \mathcal{O}_{\P^3}(-1)$,
  and $f$ contracts $D$ to a point;
\item $D \iso \P^2$ with normal bundle $\mathcal{O}_D (D) \iso
  \mathcal{O}_{\P^2}(-2)$, and $f$ contracts $D$ to a point.
\end{enumerate}
Moreover if $X$ is smooth, then in Case 1 above $Z$ is smooth,
$C$ is a smooth curve in $Z$, and $D \arrow C$ is a $\P^1$-bundle.
\end{theorem}

\begin{corollary} \label{corollary-mori} 
Suppose $X$ is a $\Q$-factorial terminal Gorenstein 3-fold and $x \in
\Mov{X}$. Then $x \cdot R \geq 0$ for any $K$-negative extremal ray
$R$ of $\Curv{X}$.
\end{corollary} {\bf Proof of Corollary:} By continuity and homogeneity it suffices
to consider the case when $x=D$ is the class of a movable Cartier
divisor. If $D \cdot R <0$ then any curve with class in $R$ must be
contained in the base locus of $D$. But by the classification in
Theorem \ref{theorem-moriclassification}, for any such $R$ the union
of all curves with class in $R$ has codimension at most 1 in $X$,
which contradicts the fact that $D$ is movable. \quad QED

Now we can prove the decomposition of the movable cone into nef cones
of SQMs. The key point is that all the SQMs $\alpha: X'
\dashrightarrow X$ we obtain are {\it over} $S$: that is, they come
with a morphism $f': X' \arrow S$ such that $f' = f \circ \alpha$. Our
proof of the theorem is very similar to that of Kawamata
\cite{Kawamata1997} for the case of a Calabi--Yau fibre space.

Before stating the theorem we need some definitions. For a normal
variety $Y$ and an $\R$-Cartier divisor $D$ on $Y$, a {\it
  $D$-flopping contraction} is a proper birational morphism $f: Y
\arrow Z$ to a normal variety $Z$ such that $\operatorname{Exc}(f)$
has codimension at least 2 in $Y$, the divisor $-(K_Y+D)$ is
$\R$-Cartier and ample over $Z$, and $K_Y$ is numerically trivial over
$Z$. The contraction $f$ is called {\it extremal} if it has relative
Picard number 1. In particular if $(Y,\Delta)$ is a $\Q$-factorial klt
pair with $\Delta$ effective and $f$ is the contraction of a
$(K_Y+\Delta)$-negative extremal ray, it is extremal, because all the
curves contracted are numerical multiples of each other, by the cone
theorem \cite[Theorem 3.7]{KollarMori1998}. Given a $D$-flopping
contraction $f: Y \arrow Z$, the {\it $D$-flop} is a birational
morphism $f^+ : Y^+ \arrow Z$ from a normal variety $Y^+$ such that
$\operatorname{Exc }(f^+)$ has codimension at least 2 in $Y^+$, and
the proper transform $(K_Y+D)^+$ of $(K_Y+D)$ on $Y^+$ is $\R$-Cartier
and ample over $Z$. (Sometimes we abuse terminology by referring to
either the birational map $f^+ \circ f^{-1} : Y \dashrightarrow Y^+$
or the variety $Y^+$ as the flop.)

Note that flops preserve the $\Q$-factorial property and for terminal
3-folds preserve the singularity type \cite[Proposition 3.37, Theorem
  6.15]{KollarMori1998}. Moreover it is easy to see that if $-K_X$ is
semiample then so is $-K_{X^+}$ for $X^+$ any flop of $X$, and these
two line bundles have the same Iitaka dimension. We conclude that the
class of 3-folds considered in Theorem \ref{theorem-maintheorem1}
(terminal $\Q$-factorial Gorenstein with $-K$ semiample of positive
Iitaka dimension) is preserved by flops.

\begin{theorem} \label{theorem-decomp} Let $X$ be a $\Q$-factorial
  terminal Gorenstein 3-fold with $-K_X$ semiample of positive Iitaka
  dimension. Then the effective movable cone $\Mov{X}^e$ decomposes as
  a union of effective nef cones of SQMs of $X$:
\begin{align*}
  \Mov{X}^e = \bigcup  \Nef{X',\alpha}^e 
\end{align*}
where the union on the right-hand side is over all SQMs $X'
\dashrightarrow X$. All these SQMs are over $S$.  The interiors of the
cones $\Nef{X',\alpha}$ are disjoint.
\end{theorem}
{\bf Proof:} Suppose $D \in \Mov{X}^e$ is an effective $\Q$-divisor on
$X$ which is not nef. By Corollary \ref{corollary-mori} $D$ cannot be
negative on a $K$-negative extremal ray, so we must have $D \cdot R
<0$ for some extremal ray $R$ of $\Curv{X}$ which lies in
$K^\perp$. Choosing some $\epsilon >0$ sufficiently small, the cone
theorem for the klt pair $(X, \epsilon D)$ tells us that $R$ is
spanned by the class of a curve, and the contraction of $R$
exists. Note that since $R \subset K^\perp$, this contraction is over
$S$. By \cite[Theorem 6.14]{KollarMori1998} the $D$-flop of this
contraction exists, and is an SQM of $X$ over $S$. If $D$ is not a
$\Q$-divisor, we can choose a small ample $R$-divisor $D'$ such that
$D+D'$ is a $\Q$-divisor in the cone $\Mov{X}^e$ but is not nef. So as
before we get a $(D+D')$-flopping contraction over $S$ and its
$(D+D')$-flop over $S$. Since $D'$ is ample, this is in particular a
$D$-flop.

So given a non-nef divisor $D \in \Mov{X}^e$, there exists a $D$-flop
$X \dashrightarrow X^+$. Applying this fact repeatedly (using the fact
explained above that $X^+$ satisfies the same assumptions as $X$),
either $D$ becomes nef after a finite sequence of $D$-flops, or else
there is an infinite sequence of $D$-flops. (Here we are abusing
terminology a little: a sequence of $D$-flops really means a sequence
of flops whose $(i+1)^{th}$ member is a $D^i$-flop, where $D^i$ is the
proper transform of $D$ by the composition of the first $i$ flops in
the sequence.) But Kawamata \cite{Kawamata1992} showed there is no
infinite sequence of $D$-flops for $D$ a $\Q$-divisor on a
$\Q$-factorial terminal 3-fold, and as remarked in \cite[Theorem
  2.3]{Kawamata1997} the same proof works if $D$ is an
$\R$-divisor. We conclude that any $D \in \Mov{X}^e$ becomes nef after
a finite sequence of flops.

We have shown that any effective movable divisor belongs to one of the
effective nef cones $\Nef{X',\alpha}^e$ where $\alpha: X'
\dashrightarrow X$ is a sequence of flops over $S$. So we have the
inclusion $\Mov{X}^e \subset \bigcup_\alpha \Nef{X',\alpha}^e$. The reverse
inclusion is clear, since an ample divisor on any SQM $X'$ is movable
on $X$, so taking closures and intersecting with the effective cone we
get $\bigcup_\alpha \Nef{X',\alpha}^e \subset \Mov{X}^e $.

To see that these flops give all the SQMs of $X$ up to isomorphism,
suppose that $\beta: Y \dashrightarrow X$ is any SQM. By the argument
above we have $\Nef{Y,\beta} \subset \bigcup_\alpha \Nef{X',\alpha}$,
so the ample cone of $Y$ must intersect the ample cone of one of the
flops, say $\alpha_i: X_i \dashrightarrow X$. So there exists a
divisor $D$ on $X$ such that ${\alpha_i}^{-1}_* D$ and $\beta^{-1}_*
D$ are ample on $X_i$ and $Y$ respectively. Therefore

\begin{align*}
X_i = \operatorname{Proj} R(X_i,{\alpha_i}^{-1}_* D) \iso
\operatorname{Proj} R(Y,\beta^{-1}_*
D) = Y
\end{align*}
and the isomorphism is compatible with $\alpha_i$ and $\beta$. In
other words, $(X_i,\alpha_i)$ and $(Y,\beta)$ are isomorphic as SQMs
of $X$.

Finally, the same argument applied to $2$ SQMs $\alpha_1: X_1
\dashrightarrow X$ and $\alpha_2: X_2 \dashrightarrow X$ shows that
the interiors of the cones $\Nef{X_1,\alpha_1}$ and
$\Nef{X_2,\alpha_2}$ are disjoint in $\Mov{X}^e$. \quad QED

As an immediate consequence, we get Theorem
\ref{theorem-maintheorem1}:

\begin{corollary}
Let $X$ be a $\Q$-factorial terminal Gorenstein 3-fold, and assume
that $-K_X$ is semiample of positive Iitaka dimension. Then the number
of effective nef cones in the decomposition $\Mov{X}^e = \bigcup_\alpha
\Nef{X',\alpha}^e$ is finite up to the action of $\PsAut(X/S)$.
\end{corollary}
{\bf Proof:} Under the hypotheses, Kawamata \cite[Theorem 3.6, Theorem
  4.4]{Kawamata1997} proves that the decomposition of $\Mov{X/S}^e$
into relative nef cones is finite up to the action of
$\PsAut(X/S)$. In other words, there are only finitely many SQMs of
$X$ over $S$ up to the action of $\PsAut(X/S)$. But Theorem
\ref{theorem-decomp} shows that every SQM of $X$ is over $S$, so we
get the result. \quad QED

\section{Extension of Kawamata's theorems}

In this section we adapt to our situation some results of Kawamata
\cite[Theorem 1.9, Theorem 2.6]{Kawamata1997} concerning the local
structure of the nef and movable cone. These results will be used in
the next section to give a description of the relative movable cone in
some cases. The proofs here are appropriately modified versions of
Kawamata's, though we give more details.

\begin{theorem} \label{theorem-localratpoly}
Let $(X/Z,\Delta)$ be a klt Calabi--Yau pair (of any dimension). Then
the cone
\begin{align*}
\Nef{X/Z} \cap B(X/Z) = \Nef{X/Z}^e \cap B(X/Z)
\end{align*}
is locally rational polyhedral inside the big cone
$B(X/Z)$. Moreover, any face $F$ of this cone corresponds to a
birational contraction $\phi: X \arrow Y$ over $Z$ by the equality $F =
\phi^*(\Nef{Y/Z} \cap B(Y/Z))$.
\end{theorem}
We will only need the theorem in the case $Z = \operatorname{Spec }k$, but
since the proof of the more general statement is identical, it makes
sense to include it here.

\noindent {\bf Proof:} The proof of the first statement works just as Kawamata's
proof for Calabi--Yau fibre spaces \cite[Theorem 5.7]{Kawamata1988},
replacing $K_X$ by $K_X+\Delta$ wherever necessary. Here is the
argument in more detail. Let $(X/Z,\Delta)$ be a klt Calabi--Yau pair,
and $D$ a $\Q$-divisor which is effective but not nef over $Z$. For a
sufficiently small positive number $\epsilon$, the pair
$(X,\Delta+\epsilon D)$ is again klt, so by the relative cone theorem
we have
\begin{align*}
\Curv{X/Z} = \Curv{X/Z} \cap \{K_X+\Delta + \epsilon D \geq 0\} + \Sigma_j R_j
\end{align*}
where the $R_j$ are extremal rays which are negative with respect to
$K_X+\Delta+\epsilon D \equiv_Z \epsilon D$. Now suppose that $D$ is
big over $Z$: then we can write $D \equiv_Z D_1 + D_2$, where $D_1$ is
effective over $Z$ and $D_2$ is ample over $Z$. Then the relative cone
theorem for the klt pair $(X,\Delta + \epsilon D_1)$ says that the
number of extremal rays $R'_j$ with $(\epsilon D_1 + \epsilon D_2)
\cdot R <0$ is finite, because $\epsilon D_2$ is ample over $Z$. But
now $(\epsilon D_1 + \epsilon D_2) \cdot R <0$ implies that $D_1 \cdot
R < 0$, so the rays $R_j$ appearing in the displayed equation above
are a subset of the set of rays $R'_j$. So there are finitely many
extremal rays of $\Curv{X/Z}$ which are negative with respect to $D$.

The nef cone of $X$ over $Z$ is the intersection $\bigcap_R \{ x \in
N^1(X/Z) | x\cdot R \geq 0\}$ of the nonnegative half-spaces in
$N^1(X/Z)$ for all extremal rays $R$ of $\Curv{X/Z}$. By the previous
paragraph, $D$ lies in all but finitely many of these
half-spaces. This means that the part of the boundary of $\Nef{X/Z}$
visible from $D$ is defined by a finite set of rational
hyperplanes. Since this is true for an arbitrary non-nef rational
point $D$ in the big cone, we conclude that $\Nef{X/Z}$ is locally
rational polyhedral inside the nef cone.

For the second statement, it suffices to show that a nef and big
divisor over $Z$ is semiample over $Z$. This will follow from the
relative version of the Basepoint-free theorem \cite[Theorem
  3.3]{KollarMori1998} if we show that $D$ being nef and big over $Z$
implies that $D-K_X-\Delta$ is nef and big over $Z$. But by definition
of a Calabi--Yau pair $K_X+\Delta = 0$ in $N^1(X/Z)$, and the property
of being nef and big is defined on the level of numerical classes
\cite[Proposition 2.61]{KollarMori1998}, so we are done. \quad QED

\begin{theorem} \label{theorem-localfiniteness}
Let $X$ be a $\Q$-factorial terminal Gorenstein 3-fold with $-K_X$
semiample of positive Iitaka dimension. Then the decomposition
\begin{align*}
  \Mov{X}^e \cap B(X) = \bigcup  \Nef{X',\alpha}^e \cap B(X)
\end{align*}
is locally finite inside the big cone $B(X)$ in the following sense:
if $\Sigma$ is a closed convex cone contained in $B(X) \cup \{0\}$,
then only a finite number of the cones $\Nef{X',\alpha}^e \cap B(X)$
intersect $\Sigma$.
\end{theorem}
{\bf Proof:} First recall from Section \ref{section-conj} that under
the stated hypotheses, there exists a $\Q$-divisor $\Delta$ on $X$
such that $(X,\Delta)$ is a klt Calabi--Yau pair. So Theorem
\ref{theorem-localratpoly} tells us that each nef cone in the
decomposition is locally rational polyhedral inside the big cone
$B(X)$.

Now let $x \in \Mov{X}^e \cap B(X)$. By Theorem \ref{theorem-decomp}
there exists an SQM $(X_0,\alpha_0)$ of $X$ such that $x \in
\Nef{X_0,\alpha_0}^e$. Let $F$ be the face of $\Nef{X_0,\alpha_0}^e$
whose interior contains $x$. I make two claims at this point. The first
claim is that the number of effective nef cones $\Nef{X_i,\alpha_i}^e$
which contain $F$ as a face is finite. The second claim is that one
can choose a small open cone $\Sigma_x$ containing $x$ such that
$\Mov{X}^e \cap \Sigma_x$ intersects only finitely many effective nef
cones $\Nef{X_i,\alpha_i}^e$. Given these two claims the result then
follows: if $\Sigma$ is a closed convex cone as in the statement of
the theorem, the cones $\{ \Mov{X}^e \cap \Sigma_x \cap \Sigma | x \in
\Sigma \}$ give an open cover of $\Mov{X}^e \cap \Sigma$. Furthermore,
since $B(X)$ is contained in the effective cone and $\Sigma$ is
contained in $B(X)$, we have the following equalities of cones:

\begin{align*}
\Mov{X}^e \cap \Sigma = \Mov{X}^e \cap B(X) \cap \Sigma = \Mov{X} \cap
B(X) \cap \Sigma = \Mov{X} \cap \Sigma.
\end{align*}
In particular $\Mov{X}^e \cap \Sigma$ is a closed cone in $N^1(X)$, so
its projectivisation is compact. Therefore the open cover we mentioned
has a finite subcover. Since each $\Sigma_x$ intersects only finitely
many nef cones, this gives the result.

It remains to prove the two claims. For the first, suppose $F$ is a face
of some nef cone $\Nef{X_0,\alpha_0}^e$. Then we have by Theorem
\ref{theorem-localratpoly} a birational morphism $f_0 : X_0 \arrow Y$
such that $F = f_0^* \, \Nef{Y}$.

I claim that the number of nef cones $\Nef{X_i/Y}^e$ for the SQMs of
$X_0$ over $Y$ is finite. To see this, choose an $f_0$-ample Cartier
divisor $D$ on $X_0$: by \cite[Lemma 6.28]{KollarMori1998} there is an
effective divisor $E$ such that $-D$ is $f_0$-linearly equivalent to
$E$. Also as explained before, since $-K_X$ is semiample we can find
an effective divisor $\Delta$ such that $(X_0,\Delta)$ is a klt
Calabi--Yau pair. So by the relative cone theorem for the klt pair
$(X_0,\Delta + \epsilon E)$ (where $\epsilon$ is a sufficiently small
positive number) the relative cone of curves $\Curv{X_0 / Y}$ is
rational polyhedral.

Now suppose there were infinitely many nef cones
$\Nef{X_i/Y}^e$. Since each $X_i$ is $\Q$-factorial and terminal, and
for each $i$ there is a divisor $\Delta_i$ such that $K_{X_i} +
\Delta_i \equiv 0$, we can apply Kawamata's decomposition theorem
\cite[Theorem 6.38]{KollarMori1998} to conclude that the birational
map $X_i \dashrightarrow X_j$ (given by the marking) between any two
SQMs decomposes as a finite sequence of flops over $Y$. (More
precisely, we are extending the result of Kawamata from the case where
$K_{X_i}$ is nef over $Y$ to the case where $K_{X_i} + \Delta_i$ is
nef over $Y$. But examining the proof we see as before that it depends
on applying the cone theorem for a klt pair $(X,\epsilon D)$, and
everything works just as well with the klt pair $(X,\Delta+\epsilon D)$.)

Now fix a model $f_0 : X_0 \arrow Y$ and an ample effective divisor
$D_0$ on $X_0$. As in the previous paragraph, there is an effective
divisor $E_0$ on $X_0$ which is $f_0$-linearly equivalent to
$-D_0$. Then any other $X_i$ is obtained from $X_0$ by a finite
sequence of $E_0$-flops. Since each cone $\Curv{X_i / Y}$ is rational
polyhedral, there are only finitely many ways to do a flop at each
stage. But we assumed that there are infinitely many nef cones, so by
K\"onig's Lemma, there must be an infinite sequence of
$E_0$-flops. This contradicts the fact mentioned in the proof of
Theorem \ref{theorem-decomp} that any sequence of terminal 3-flops is
finite. So there can be only finitely many nef cones
$\Nef{X_i/Y}^e$. Call the corresponding SQMs $X_1, \ldots, X_n$.

Finally suppose that $X_\alpha$ is any SQM of $X$ such that
$\Nef{X_\alpha}^e$ contains $F$ as a face. As before, we get a
birational contraction $f_\alpha: X_\alpha \arrow Y$. By the previous
paragraph, we must have $\Nef{X_\alpha/Y}^e = \Nef{X_i/Y}^e$ for some
$i \in \{1,\ldots,n\}$. This means $X_\alpha$ is isomorphic to $X_i$
as an SQM over $Y$, so in particular isomorphic to $X_i$ as an SQM of
$X$ (over $\operatorname{Spec} k$). In other words $\Nef{X_\alpha}^e =
\Nef{X_i}^e$ as cones in $N^1(X)$, which completes the proof of the
first claim.

Now for the second claim: namely, given $x \in \Mov{X}^e \cap B(X)$
one can choose a small open cone $\Sigma_x$ containing $x$ such that
$\Mov{X}^e \cap \Sigma_x$ intersects only finitely many effective nef
cones $\Nef{X_i}^e$. Recall that $x$ belongs to the interior of a face
$F$ of some nef cone, and we have just shown that there are only
finitely many nef cones which contain $F$ as a face. Call these nef
cones $\Nef{X_i}$, for $i = 1,\ldots,k$. By Theorem
\ref{theorem-localratpoly} each nef cone is locally rational polyhedral
near $x$, so for any sufficiently small open cone $K$ containing $x$,
the only faces of these nef cones which meet $K$ are the faces which
contain $F$ as a face. 

If the claim is not true, then any open neighbourhood of $x$ inside
$\Mov{X}$ must meet infinitely many nef cones, and in particular must
meet infinitely many nef cones which do not contain $F$ as a face
(because there are only finitely many nef cones which do contain $F$
as a face). Let us assume this is the case, and derive a
contradiction.

Let $K$ be an open cone around $x$ chosen small enough (as explained
above) so that the only faces of the nef cones $\Nef{X_i}$ (for
$i=1,\ldots,k$) which meet $K$ are those which contain $F$ as a
face. I claim that the union $U = \Nef{X_1} \cup \cdots \cup \Nef{X_k}$
must have an `exposed face' inside $\Mov{X} \cap K$: that is, there
is a face $F^\prime$ of one of these nef cones which lies on the
relative boundary inside $\Mov{X} \cap K$ of the set $U$. To see this,
note that if no such exposed face exists, then every point of $U \cap
K$ would be in the relative interior of $U \cap K$ inside $\Mov{X}
\cap K$, implying that $U \cap K$ is an open subset of $\Mov{X} \cap
K$. On the other hand , $U \cap K$ is a finite union of closed sets in
$\Mov{X} \cap K$, so is closed. Therefore $U \cap K$ must be the whole
of $\Mov{X} \cap K$. So $\Mov{X} \cap K$ is an open neighbourhood of
$x$ inside $\Mov{X}$ which meets only finitely many nef cones,
contradicting the hypothesis of the previous paragraph. This
contradiction implies that $U$ must have an exposed face inside
$\Mov{X} \cap K$.

Since the relative boundary of $U$ inside $\Mov{X} \cap K$ is itself a
polyhedral complex, there must in fact be an exposed face of
codimension 1 inside $\Mov{X} \cap K$. Choose such a face $F'$; by our
choice of the cone $K$ in the previous paragraph we know that $F'$
contains $F$ as a face. Assume that $F'$ is a face of a nef cone
$\Nef{X_1}$. Then $F'$ is dual to an extremal ray $R$ of the closed
cone of curves $\Curv{X_1}$. Choose an interior point $x'$ of $F'$:
since $F'$ is exposed, any open neighbourhood of $x'$ inside $\Mov{X}
\cap K$ contains points which do not belong to the set $U$. Choose a
point $x''$ in the complement of $U$ with the property that $x''$ is
negative on the extremal ray $R$, but is positive on all other
extremal rays of $\Curv{X_1}$. (Here we are using the conclusion of
Theorem \ref{theorem-localratpoly} that $\Nef{X_1}$ is locally
rational polyhedral inside the big cone, so in a neighbourhood of $x'$
it is defined by finitely many hyperplanes. By the standard argument
with the cone theorem (applied to a klt pair $(X_1,\epsilon x'')$ for
$\epsilon$ sufficiently small), we obtain the contraction of the
extremal ray $R$, and the corresponding flop. Call the flopped variety
$X'$: then $\Nef{X'}$ intersects $\Nef{X_1}$ along the face $F'$, so
in particular contains $F$ as a face. In other words, $\Nef{X'}$ must
equal $\Nef{X_i}$ for some $i=1,\ldots,k$.

On the other hand since the two nef cones $\Nef{X'}$ and $\Nef{X_i}$
meet along the codimension-1 face $F'$, the union $\Nef{X_1} \cup
\Nef{X'}$, which we now know is a subset of $U$, contains an open
neighbourhood of every interior point of $F'$. This contradicts the
fact that $F'$ is an exposed face of $U$.

The contradiction shows that our original assumption that every open
neighbourhood of $x$ inside $\Mov{X}$ meets infinitely many nef cones
is untenable, and so the second claim is proven. \quad QED

\section{Relative movable cone and lifting} \label{section-movablecone}

In this section we describe the relative movable cone completely in
the case that $-K_X$ has Iitaka dimension 2 and the fibration $f: X
\arrow S$ has only 1-dimensional fibres. Roughly speaking, we show
that the relative movable cone is `as simple as possible': it is
bounded by the numerical classes of moving families of fibre
components. The precise statement is given in Theorem
\ref{theorem-relativemovable}. This description allows us to prove a
finiteness result for the (absolute) movable cone: there is a rational
polyhedral cone which intersects only finitely many nef cones, and
whose translates cover the whole movable cone (Theorem
\ref{theorem-cone}). As an application, we get a proof of the full
movable cone conjecture in a special case (Example
\ref{example-quadrics}). Throughout this section we work under the
following assumption:

\vspace{11pt}
\noindent {\bf Assumption 1:} {\it $X$ is a $\Q$-factorial terminal
  Gorenstein 3-fold, the anticanonical bundle $-K_X$ is semiample of
  Iitaka dimension 2, and the fibres of the contraction morphism
  $f:~X~\arrow~S$ defined by $-K_X$ are all 1-dimensional.}

\vspace{11pt}
\noindent It is not clear to me to what extent the methods of this
section can be adapted to study the non-equidimensional case, or the
case when $-K_X$ has Iitaka dimension 1. In any case, the conditions
of Assumption 1 still permit at least one nontrivial class of
examples, discussed in Example \ref{example-quadrics}.

We begin by introducing some notation and terminology. We denote the
generic fibre of the morphism $f: X \arrow S$ by $X_\eta$. Since $f$
is given by some sections of some power of $-K_X$, adjunction tells us
that $X_\eta$ is a smooth curve of genus 1, and so we refer to $f$ as
an {\it elliptic fibration}. (Note that we do not assume the existence
of a rational section of $f$.)  The {\it Mordell--Weil group} of $f$
is the abelian group $\Pic^0(X_\eta)$ of degree-0 line bundles on
$X_\eta$. The translation action of the Mordell--Weil group on
$X_\eta$ extends to an action on $X$ by pseudo-automorphisms over $S$,
so we have an inclusion $\Pic^0(X_\eta) \subseteq \PsAut(X/S)$. To
prove any form of the movable cone conjecture, it is therefore
sufficient to prove the corresponding statement with $\Pic^0(X_\eta)$
in place of $\PsAut(X/S)$, and this is what Kawamata does in his proof
for Calabi--Yau fibre spaces.

We start by quoting the following lemma \cite[Lemma 3.1]{Kawamata1997}
which allows us to pull back arbitrary Weil divisors from $S$ to $X$.

\begin{lemma} \label{lemma-qfacbase}
Let $f: X \arrow S$ be as in Assumption 1. Then $S$ is
$\Q$-factorial. 
\end{lemma}

Next we define the subspace $T(X/S)$ of {\it relatively trivial}
divisors on $X$ to be the kernel of the natural projection $N^1(X)
\arrow N^1(X/S)$. In other words $T(X/S)$ is the subspace of classes
which have degree zero on all curves in $X$ which map to a point in
$S$. It is clear that $T(X/S)$ contains $f^*(N^1(S))$ as a subspace,
but in general they need not be equal \cite[Proposition
  3.27]{KollarMori1998}. The first thing to note is that $T(X/S)$ does
not depend on the model of $X$ we choose:

\begin{lemma} \label{lemma-numtrivial}
    Suppose $X' \arrow S$ is any SQM of $X$ (automatically over $S$ by
    Theorem \ref{theorem-decomp}). Then $T(X/S)=T(X'/S)$ as subspaces
    of $N^1(X) = N^1(X')$.
\end{lemma}
{\bf Proof:} By Theorem \ref{theorem-decomp} $X'$ is obtained from $X$
by flopping fibral curves. Suppose first $X \dashrightarrow X'$ is the
flop of a single curve $C$ with numerical class $[C]$. The canonical
identification $N^1(X) \iso N^1(X')$ gives a dual identification
$N_1(X) \iso N^1(X')$ which takes $[C]$ to $-[C]$ and, for any other
curve $\Gamma$ in $X$, takes $[\Gamma]$ to $[\Gamma]+n[C]$ for some
integer $n$ (depending on $\Gamma$). Therefore the identification
takes the subspace $N_1(X/S)$ to $N_1(X'/S)$.

If $X \dashrightarrow X'$ is any sequence of flops, applying the same
argument repeatedly shows that $N_1(X/S)$ is identified with
$N_1(X'/S)$.  Now $T(X/S) = N_1(X/S)^\perp \subset N^1(X)$, giving the
result. \quad QED

Next, a prime divisor $D$ on $X$ is defined to be {\it vertical} if
$f(D) \neq S$. We denote by $V(X/S)$ the subspace of $N^1(X/S)$
spanned by vertical divisors, and denote by $v$ the dimension of
$V(X/S)$. A prime divisor $D$ is {\it exceptional over $S$} if there
is a SQM $f':X' \arrow S$ of $X$ over $S$ and a divisorial contraction
$\phi: X' \arrow Y$ over $S$ with exceptional divisor $E$ such that
$D$ is the proper transform of $E$. We will need the following facts
about vertical and exceptional divisors \cite[Lemma 3.1, Lemma
  3.2]{Kawamata1997}:

\begin{lemma} \label{lemma-kawamatafacts} Let $f:~X~\arrow~S$ be as in
  Assumption 1. 

1. If $D$ is an vertical prime divisor, and $D$ does not contain the
fibre of the restriction $f^{-1}f(D) \arrow D$ over the generic point
of the curve $f(D)$, then $D$ is exceptional over $S$.

2. The vector space $V(X/S)$ of vertical divisor classes is spanned by
the classes of exceptional divisors over $S$.
\end{lemma}
Since we assume that all fibres of $f: X \arrow S$ are 1-dimensional,
all fibres have the same numerical class in $N_1(X)$, which we denote
by $F$. Now suppose $C$ is an irreducible curve on $S$. The preimage
$f^{-1}(C)$ of this curve is (as a set) a union $D_1 \cup D_2 \cup
\cdots \cup D_k$ of vertical prime divisors on $X$. The set-theoretic
intersection $D_i \cap f^{-1}(p)$ of the divisor $D_i$ with a fibre of
$f$ over the generic point of $C$ is a curve whose numerical class we
denote by $F_i$. Since the numerical class of any fibre is equal to
$F$, we get an expression $m_1F_1 + \cdots + m_kF_k = F$, where the
$m_i$ are positive integers. This implies that $D_i \cdot (m_1F_1 +
\cdots + m_kF_k) =0$ for each $i$, because $D_i$ is vertical. On the
other hand since the fibres of $f$ are connected, if $k \geq 2$ (which
by Lemma \ref{lemma-kawamatafacts} means that the $D_i$ are
exceptional over $S$) we must have $D_i \cdot F_j >0$ for some $j \neq
i$, and therefore $D_i \cdot F_i <0$. This shows in particular that
the classes $F_i$ for $i=1,\ldots, k$ are all distinct in
$N_1(X)$. Furthermore if $D_1$ and $D_2$ are exceptional prime
divisors over $S$ such that $f(D_1) \neq f(D_2)$, with the
corresponding curves having classes $F_1$ and $F_2$ then we have $D_1
\cdot F_2 = D_2 \cdot F_1 =0$, again implying that $F_1$ and $F_2$ are
distinct classes. We conclude that for each exceptional prime divisor
$D_i$ over $S$ we have a class $F_i \in N_1(X)$ of a fibral curve of
$f$, and distinct divisors give distinct classes. Since $f$ has
1-dimensional fibres, there can be only finitely many such classes
(because there are only finitely many decompositions of $F$ in the
monoid $\text{Curv}(X)_\Z$), so there are only finitely many
exceptional divisors over $S$. From now on $D_1,\ldots,D_n$ will
denote these exceptional divisors over $S$, and $F_1,\ldots,F_n$ the
corresponding classes of fibral curves of $f$.

One other piece of notation will be useful for us. We partition the
set $\mathcal{D} = \{D_1,\ldots,D_n\}$ of exceptional divisors over
$S$ into subsets $\mathcal{D}_1,\ldots,\mathcal{D}_r$ by saying that
$D_i$ and $D_j$ belong to the same member of the partition if the
curves $f(D_i)$ and $f(D_j)$ are the same. Denote by $C_p$ the curve
in $S$ which is the common image of the divisors belonging to the
member $\mathcal{D}_p$ of the partition. Note that $f^*(C_p)=
\sum_{\{i \, | \, D_i \in \mathcal{D}_p\}} \mu_i D_i$ for some
positive rational numbers $\mu_i$, so in particular $\sum_{\{i \, | \,
  D_i \in \mathcal{D}_p\}} \mu_i D_i=0$ in $N^1(X/S)$. Also, by the
previous paragraph, for each $\mathcal{D}_p$ there are positive
integers $m_i$ such that $\sum_{\{i \, | \, D_i \in \mathcal{D}_p\}}
m_i F_i =F$, the class of a fibre. 

Now we can begin to describe the relative movable cone of $X$ over
$S$. We need the following `negativity lemma' for sums of exceptional
divisors over $S$.

\begin{lemma} \label{lemma-hodgeindex}
Suppose $D_i$ and $F_j$ are as defined above. Then a class $x = \sum_i
r_i D_i \in N^1(X)$ satisfies $x \cdot F_j \geq 0$ for all
$j=1,\ldots,n$ if and only if $x = f^*(C)$ for some $\R$-divisor $C$
on $S$, in which case $x \cdot F_j=0 $ for all $j$. A class $y =
\sum_i s_i F_i \in N_1(X)$ satisfies $D_j \cdot y \geq 0$ for all
$j=1,\ldots,n$ if and only if $y =sF$ for some real number $s$, in
which case $D_j \cdot y =0$ for all $j$.
\end{lemma}
{\bf Proof:} One direction is trivial: if $x$ is the pullback of an
$\R$-divisor on $S$ the intersection numbers $x \cdot F_j$ are
obviously 0. Conversely assume $x = \sum_i r_i D_i$ satisfies $x \cdot
F_j \geq 0$ for all $j=1,\ldots,n$. We want to show that $x$ is the
pullback of an $\R$-divisor on $S$.

First, for each $\mathcal{D}_p$ in the partition defined above there
are positive integers $m_i$ such that $\sum_{\{i \, | \, D_i \in
  \mathcal{D}_p\}} m_i F_i =F$. Therefore ince $D_i \cdot F =0$ for
all $D_i$, the condition $x \cdot F_j \geq 0$ for all $j$ is
equivalent to the condition $x \cdot F_j = 0$ for all $j$.

Next, the partition we defined gives a decomposition of $x$ as
$x=x_1+\cdots+x_r$, where $x_p = \sum_{\{i \, | \, D_i \in
  \mathcal{D}_p\}} r_i D_i$. Note that $x_p \cdot F_j =0$ whenever $j$
is not in $\mathcal{D}_p$, so the condition $x \cdot F_j = 0$ for each
$j$ is equivalent to the following: for each $j$, we have $x_p \cdot
F_j = 0$ for the unique $p$ such that $j \in \mathcal{D}_p$. So we can
reduce to the case where $x=x_p$ for some $p$: that is, $x$ is a real
linear combination of divisors $D_i$ all belonging to some member
$\mathcal{D}_p$ of the partition, that is all of whose images $f(D_i)
\subset S$ are the same curve. Call this curve $C$, and write $f^*(C)
= \sum_{\{i \, | \, D_i \in \mathcal{D}_p\}} \mu_i D_i$. We want to
show that there exists a real number $\lambda$ such that $r_i =
\lambda \mu_i$ for all $i$, implying $x = \lambda f^*(C)$.

To see this, choose a general curve $\Gamma \subset S$ intersecting
$C$ transversely in $m>0$ points. We will restrict to the surface
$\Sigma = f^{-1}(\Gamma)$. The restriction of $D_i$ to $\Sigma$ is a
union of $m$ disjoint curves each with class $F_i$, so we get
$x_{|\Sigma} = m(\sum_i r_i F_i)$. Since $x \cdot F_j=0$ for each $j$,
we get $(x_{|\Sigma})^2 =0$. But the Hodge index theorem on $\Sigma$
then implies that $x_{|\Sigma}$ is a multiple of the class $\Phi$ of a
fibre of $\Sigma \arrow \Gamma$: say $x_{|\Sigma} = \lambda \Phi$. By
looking at the fibre over a point of $C$ we see that $\Phi = \sum_i
\mu_i F_i$, so we conclude that $\sum_i r_i F_i = \lambda (\sum_i
\mu_i F_i)$, hence $x = \lambda f^*(C)$, as required.

The proof of the claim about the classes $y = \sum_i s_i F_i$ is analogous.
\quad QED

\begin{corollary} \label{corollary-dualspaces}
The dual vector space
\begin{align*} F(X/S) = N_1(X/S)/ \{ x \in N_1(X/S) | D \cdot
x =0 \text{ for all } D \in V(X/S) \}
\end{align*}
 of the space $V(X/S)$ of vertical divisors is spanned by the classes
 $F_1,\ldots,F_n$.
\end{corollary}

\begin{theorem} \label{theorem-relativemovable}
Let $D_i$ and $F_j$ be as defined above. Then the relative effective
cone and the relative movable cone of $X$ over $S$ are as follows:

\begin{align*}
  B^e(X/S) &=\{x \in N^1(X/\P^2) : x \cdot F >0 \} \cup
  \R_+\{[D_i]\} \cup \{0\} \\ \Mov{X/\P^2}^e &= \{x \in
  N^1(X/S) : x \cdot F >0, \, x \cdot F_i \geq 0 \,
  (i=1,\ldots,n) \} \cup \{0\}.
\end{align*}
\end{theorem}
{\bf Proof:} First suppose that $D$ is an effective class over $S$
with degree $k \leq 0$ on the generic fibre.  There exists a nonempty
open set $U \subset S$ such that $D(f^{-1}(U)) \neq 0$. Choose a
nonzero section $s \in D(f^{-1}(U))$. Then the class of the divisor
$\Delta=\overline{\{s=0\}}$ differs from $D$ only on the codimension-1
subset $X \backslash f^{-1}(U)$: that is, there is an effective
divisor $D'$ supported on $X \backslash f^{-1}(U)$ such that $\Delta =
D+D'$ in $N^1(X)$. In particular $D$ and $\Delta$ have the same degree
$k$ on the generic fibre. Since $\Delta$ is effective this implies
$k=0$.  Moreover $k=0$ implies that $\Delta$ is a sum of vertical
(prime) divisors so by Lemma \ref{lemma-kawamatafacts} its class
belongs to the cone $V$ spanned by $f^*(B^e(S))$ and the
$D_i$. Finally the divisor $\Delta-D$ is supported in $X \backslash
f^{-1}(U)$ therefore its support maps onto a curve in $S$. So any
divisor in the support of $\Delta-D$ must also have class in the cone
$V$. So for any effective class $D$ over $S$ with degree $\leq 0$ on
the generic fibre we can write $D=V_1-V_2$ where $V_i$ are classes in
$V$. The image of $V$ in $N^1(X/S)$ is the cone $\R_+\{[D_i]\}$, which
is closed under negation since for any $i$, there is member
$\mathcal{D}_p$ of the partition which contains $D_i$, and there exist
positive rational numbers $m_k$ such that $\sum_{k: D_k \in
  \mathcal{D}_p} m_k D_k=0$ in $N^1(X/S)$. Therefore $[D]=[V_1]-[V_2]$
belongs to this cone as claimed. This proves that the left-hand side
of the first equation is contained in the right-hand side. To prove
the reverse inclusion, first note that if a divisor $D$ has positive
degree on an irreducible fibre $F$ then the restriction $D_{|F}$ is
ample hence effective. But standard results on semicontinuity of
cohomology \cite[Corollary III.12.9]{Hartshorne1977} show that for any
positive integer $m$, any section of $mD_{|F}$ is the restriction of a
section in $mD(f^{-1}(U))$ for $U \in S$ some open subset. Therefore
by definition $mD$ and hence $D$ is effective over $S$. Finally, all
divisors in the cone $V$ are effective by definition hence
effective over $S$, so all elements of $\R_+\{[D_i]\}$ lie in the
relative effective cone. This completes the proof of the claim about the
relative effective cone.

Now we must prove the claim about the relative movable cone. First
note that if $D$ is a class in $N^1(X)$ with $D \cdot F_i <0$ for some
$i$, then $D$ cannot be movable over $S$. For suppose $C$ is a curve
in $X$ with class $F_i$. If there was an open set $U \subset S$
containing the point $f(C)$ and a section of $D(f^{-1}(U))$ not
vanishing identically along $C$ we would have $D \cdot C \geq 0$
contradicting our assumption: therefore every such curve $C$ is
contained in $\operatorname{Supp \ Coker } (f^*f_* \mathcal{O}_X(D)
\arrow \mathcal{O}_X(D))$. Since these curves $C$ fill up an open set
in the divisor $D_i$ we conclude that $D$ cannot be movable over
$S$. So the relative movable cone is contained in the cone $\{x \cdot
F_i \geq 0 \text{ for all } i\}$. If moreover $x$ is a nonzero class
effective and movable over $S$, it must have $x \cdot F>0$. For
otherwise by the description of $B^e(X/\P^2)$ we would have $x \in
\R_+\{[D_i]\}$. Any nonzero point in this cone can be written in the
form $x = \sum r_i[D_i]$; by Lemma \ref{lemma-hodgeindex}, since $x$
is not pulled back from $S$, this implies $x \cdot F_i <0$ for some
$i$. So we have shown that left-hand side is contained in the
right-hand side in the second equality above.

Conversely suppose that $x \in N^1(X/S)$ satisfies $x \cdot F_i \geq
0$ for all $i$ and $x \cdot F>0$: we want to show that $x$ belongs to
the relative effective movable cone. First note that any such $x$ is
effective over $S$ by our description of the relative effective
cone. Next suppose that $D$ is a divisor class with $D \cdot F_i >0$
for each $i$. Since the class of a fibre $F$ is a sum of classes
$F_i$, the restriction of such a $D$ to any irreducible fibre is
ample. Also since $D \cdot F_i>0$ for each $i$ the restriction of $D$
to each component of the fibre over the generic point of each curve
$f(D_i)$ is ample. Again by the semicontinuity result mentioned above,
for any integer $m>0$ any section of the line bundle $mD_{|F}$ comes
from a section in $mD(f^{-1}(U))$ for some open $U \subset S$
containing $f(F)$. Therefore if we choose an integer $m$ sufficiently
large so that $mD_|F$ is very ample for some fibral curve $F$, then
$\operatorname{Supp \ Coker } (f^*f_* \mathcal{O}_X(mD) \arrow
\mathcal{O}_X(mD))$ does not contain any point in the curve $F$. We
conclude that the relative stable base locus $Bs_\Q(D,f) := \cap_{m
  \geq 1} \operatorname{Supp \ Coker } (f^*f_* \mathcal{O}_X(mD)
\arrow \mathcal{O}_X(mD))$ does not contain any point in any
irreducible fibre or any point of the fibre over the generic point of
one of the curves $f(D_i)$.

Suppose that the set $Bs_\Q(D,f)$ contains a prime divisor
$\Delta$. We have just seen that $\Delta$ must be disjoint from all
irreducible fibres of $f$, so it must be vertical. Also $\Delta$ is
disjoint from all components of the fibre of $f$ over the generic
point of $f(D_i)$, so it is not one of the $D_i$. Since the fibre of
$f$ over every point of $f(\Delta)$ is reducible, I claim that $\Delta
\neq f^{-1}f(\Delta)$, as sets. To prove this, suppose equality holds:
then $\Delta$ would contain the fibre of $f$ over every point of
$f(\Delta)$. Suppose the fibre of $f$ over the generic point of
$f(\Delta)$ is a union $C_1 \cup \cdots \cup C_r$ of irreducible
curves (with $r>1$). Taking the closure of the curve $C_i$ gives a
divisor $\Delta_i$ properly contained in $\Delta$. But by assumption
$\Delta$ is prime, a contradiction. So by Kawamata \cite[Lemma
  3.1]{Kawamata1997} the divisor $\Delta$ must be exceptional over
$S$, therefore equal to one of the $D_i$. This is again a
contradiction. We conclude that no such divisor $\Delta$ can exist, so
$Bs_\Q(D,f)$ contains no divisors. By the relative version of
`stabilisation of the base locus' \cite[Proposition
  2.1.21]{Lazarsfeld2004} $D$ is therefore movable over $S$.

We have shown that any Cartier divisor class $[D]$ with $[D] \cdot F
>0$ and $[D] \cdot F_i >0$ for all $i$ belongs to the relative movable
cone, so the same is true for rational classes. To complete the proof
we observe that any point in the cone $\{x \in N^1(X/\P^2) : x \cdot F
>0, \, x \cdot F_i \geq 0 \text{ for all } i\}$ is the limit of
rational classes $x_\alpha$ with $x_\alpha \cdot F_i >0$ for all
$i$. We have just proved that each class $x_\alpha$ belongs to the
closed cone $\Mov{X/S}$ and therefore so does their limit $x$. \quad
QED

\begin{lemma} \label{lemma-surjective}
The projection $\Mov{X}^e \arrow \Mov{X/S}^e$ is surjective.
\end{lemma}
{\bf Proof:} 
 By Theorem \ref{theorem-relativemovable} we know that $\Mov{X/S}^e$
 is generated as a convex cone by the classes of Cartier divisors, so
 it is enough to prove that if $D$ is a nonzero Cartier divisor whose
 class in $N^1(X/S)$ lies in $\Mov{X/S}^e$, then $mD+\sum_k \nu_k
 f^*(A_k)$ is effective and movable for some integers $m>0$ and
 $\nu_k$ and some ample divisors $A_k$ on $S$.

By our description of $\Mov{X/S}^e$ we see that any nonzero class in
that cone is big over $S$; also, by Theorem \ref{theorem-decomp} we
know that any class in that cone belongs to $\Nef{X'/S}^e$ for some
SQM $X'$ of $X$. The statement of the lemma is unaffected if we
replace $X$ by some SQM, so we can assume without loss of generality
that $D$ is nef over $S$ and big over $S$. By the relative
Basepoint-Free theorem, $D$ is then semiample over $S$: that is, for
some positive integer $m$ the morphism $f^*f_* \curlyo(mD) \arrow
\curlyo (mD)$ is a surjection. In particular, we can choose an open
set $U \subset S$ such that the restriction of $mD$ to $f^{-1}(U)$ is
basepoint-free. Therefore we can choose 2 sections $s_1, \, s_2 \in
mD(f^{-1}(U))$ such that the common vanishing set of $s_1$ and $s_2$
has codimension 2 in $f^{-1}(U)$. Now suppose these sections have
poles of order $m_1$ and $m_2$ along the divisor $\Delta=X \backslash
f^{-1}(U)$. Note that $\Delta$ pulls back from a divisor on $S$, so
its class belongs to the subspace $f^*(N^1(S)) \subset N^1(X)$. Since
$N^1(S)$ is spanned by the classes of ample divisors, we can write
(modulo numerically trivial divisors) $\Delta= \sum_k n_k f^*(A_k)$
where the $n_k$ are integers (not necessarily positive) and the $A_k$
are ample divisors on $S$. Now multiplying $s_1$ and $s_2$ by
appropriate powers of a section $\sigma$ of a line bundle numerically
equivalent to $\sum_k n_k f^*(A_k)$ which satisfies $\sigma^{-1}(0) =
\Delta$ (as sets), we get global sections $\sigma_1$ and $\sigma_2$ of
line bundles numerically equivalent to $mD+m_1(\sum_k n_k f^*(A_k))$
and $mD+m_2(\sum_k n_k f^*(A_k))$, and these sections do not vanish
along $\Delta$. Moreover, the restrictions of $\sigma_1$ and
$\sigma_2$ to $f^{-1}(U)$ have the same vanishing sets as $s_1$ and
$s_2$, so we conclude that the common vanishing set of $\sigma_1$ and
$\sigma_2$ has codimension 2 in $X$. Finally choose a positive
integer $n$ sufficiently large so that $n-m_1n_k$ and $n-m_2n_k$ are
positive for all $k$, and choose sections $\tau_1$ and $\tau_2$ of the
ample line bundles $\sum_k(n-m_1n_k)f^*(A_k)$ and
$\sum_k(n-m_1n_k)f^*(A_k)$ whose common vanishing set has codimension
at least 2 and such that the common vanishing set of $\tau_i$ and
$\sigma_j$ has codimension 2 for $i \neq j$. Then $\sigma_1\tau_1$ and
$\sigma_2\tau_2$ are sections of a line bundle numerically equivalent
to $mD+n(\sum_k f^*(A_k))$ whose common vanishing set has codimension
at least $2$. Therefore the class $mD+n(\sum_k f^*(A_k))$ is movable.

Now suppose $D$ is any Cartier divisor effective over $S$. By definition
this means there is an open set $U \subset S$ and a nonzero section
$s$ of $D(f^{-1}(U))$. Suppose that $s$ has a pole of order $m$ along
the divisor $\Delta=X \backslash f^{-1}(U)$. As above, multiplying $s$
by an appropriate section of a bundle numerically equivalent to
$\sum_k n_k f^*(A_k)$ we get a global section $\sigma$ of a bundle
numerically equivalent to $D+\sum_k n_k f^*(A_k)$. Therefore $D+\sum_k
n_k f^*(A_k)$ is effective, so the map $B^e(X) \arrow B^e(X/S)$ from
the effective cone to the relative effective cone is surjective.

Finally suppose $[D] \in \Mov{X/S}^e$ is the image of a Cartier
divisor $D$. We have shown that there are classes $D_1 = mD+n(\sum_k
f^*(A_k))\in \Mov{X}$ and $D_2 = D+\sum_k n_k f^*(A_k) \in B^e(X)$. So
put $D_3 = mD + \sum_k \nu_k f^*(A_k)$, where $\nu_k =
\operatorname{max }(n, mn_k)$. This gives 
\begin{align*}
D_1 + \sum_k (\nu_k -n) f^*(A_k) = D_3 = mD_2 +\sum_k (\nu_k -mn_k)
f^*(A_k).
\end{align*}
For any $r \geq 0$ and any ample divisor $A$ on $S$, the class
$rf^*(A)$ is in the intersection $B^e(X) \cap \Mov{X}$, so the
left-hand side of the displayed equation is an element of $\Mov{X}$
and the right-hand side is an element of $B^e(X)$. Therefore $D_3$ is
an element of the intersection $B^e(X) \cap \Mov{X} = \Mov{X}^e$. So
every class $[D] \in \Mov{X/S}^e$ which is the image of a Cartier
divisor class is in the image of the cone $\Mov{X}^e$, as
required. \quad QED

Having described the relevant cones, we analyse the action of the
pseudo-automorphism group on them. First we state a result of Kawamata
\cite[Lemma 3.5]{Kawamata1997} which shows that, passing to a suitable
quotient, the Mordell--Weil group of our elliptic fibration acts as a
group of translations. 

\begin{theorem}[Kawamata] \label{theorem-kawamata}
The group $\Pic^0(X_\eta)$ acts properly discontinuously on the affine
subquotient space $W(X/S) = \{ z \in N^1(X/S)/V(X/S) : z \cdot F=1 \}$
as a group of translations, with fundamental domain a rational
polyhedron.
\end{theorem}
We use this theorem together with our description of $\Mov{X/S}^e$ to
find a rational polyhedral cone whose translates by the Mordell--Weil
group cover the relative movable cone:

\begin{lemma} \label{lemma_funddomain}
There is a rational polyhedral subcone $K$ of $\Mov{X/S}^e$ such that
$\Pic^0(X_\eta) \cdot K = \Mov{X/S}^e$.
\end{lemma} {\bf Proof:} Let $W'(X/S)$ denote
the affine subspace $\left\{y \in N^1(X/S) : y \cdot F = 1 \right\}$
and denote by $q$ the quotient map $W'(X/S) \arrow W(X/S)$. By
definition of the quotient action of $\Pic^0(X_\eta)$, for any $\phi
\in \Pic^0(X_\eta)$ and $x \in N^1(X/S)$ we have
$\phi(q(x))=q(\phi(x))$. By Theorem \ref{theorem-kawamata} the action of
$\Pic^0(X_\eta)$ on $W(X/S)$ has fundamental domain a rational
polyhedron $\Pi$, and hence for the action on $W'(X/S)$ we have
$\Pic^0(X_\eta) \cdot q^{-1}(\Pi) = W(X/S)$. Since the action of
$\Pic^0(X_\eta)$ preserves the relative effective movable cone, we can
intersect with that cone on both sides to get $\Pic^0(X_\eta) \cdot
(q^{-1}(\Pi) \cap \Mov{X/S}^e) = \Mov{X/S}^e \cap W(X/S)$. Finally
since $\Pic^0(X_\eta)$ acts linearly we can multiply on both sides by
positive scalars to get $\Pic^0(X_\eta) \cdot \R_+ (q^{-1}(\Pi) \cap
\Mov{X/S}^e) = \Mov{X/S}^e$. (Here we use the fact that every ray of
$\Mov{X/S}^e$ intersects $W(X/S)$, which follows from Theorem
\ref{theorem-relativemovable}.) So taking $K =\R_+ (q^{-1}(\Pi) \cap
\Mov{X/S}^e)$ it remains to show that $q^{-1}(\Pi) \cap \Mov{X/S}^e$
is a rational polyhedron in $W'(X/S)$. Since $\Pi$ is a rational
polyhedron and by Theorem \ref{theorem-relativemovable} the cone
$\Mov{X/S}^e$ is defined by a finite set of inequalities, we need to
show that $q^{-1}(\Pi) \cap \Mov{X/S}^e$ is bounded. Choosing a
section $s$ of $q$ we can write $W'(X/S)=V(X/S) + \operatorname{im}
s$. Let $\Pi'$ denote the polyhedron $s(\Pi)$: then $q^{-1}(\Pi) =
V(X/S) + \Pi' \subset W'(X/S)$. So suppose a vector $v+s$ with $v \in
V(X/S)$ and $s \in \Pi'$ belongs to $q^{-1}(\Pi) \cap \Mov{X/S}^e$. By
Theorem \ref{theorem-relativemovable} the intersection numbers $(v+s)
\cdot F_i$ must be nonnegative for all $i$. Now $s \cdot F_i$ is
bounded for $s \in \Pi'$ by compactness of $\Pi'$, so $v \cdot F_i$ is
bounded below for all $i$. Now recall that for any member
$\mathcal{D}_p$ of our partition we have an expression $\sum_{j \in
  \mathcal{D}_p} m_j F_j =F$ for some positive integers
$m_j$. Therefore $v \cdot m_iF_i = - \sum_{j \in \mathcal{D}_p, j \neq
  i} v \cdot m_jF_j$ for any $v \in V(X/S)$, so $v \cdot F_i$ is
bounded above and below. Now fix a basis for $V(X/S)$ consisting of
the classes in $N^1(X/S)$ of some of the vertical divisors $D_i$:
relabelling, we can call these basis elements $[D_1], \ldots,
[D_v]$. For a vector $v = \sum_{i=1}^v a_i [D_i]$ we can solve for the
coefficients $a_i$ in terms of the intersection numbers $v \cdot F_i$,
because by Corollary \ref{corollary-dualspaces} the dual of $V(X/S)$
is spanned by the $F_i$. Therefore the coefficients of $v$ are
bounded. Since $\Pi'$ is compact, we conclude that the subset
$q^{-1}(\Pi) \cap \Mov{X/S}^e$ is bounded, hence rational polyhedral,
as required. \quad QED

Now we complete the proof of our main finiteness result for the
movable cone. The final step is to lift the cone $K$ from the previous
lemma to a cone in $\Mov{X}^e$ which has the properties we want.

\begin{theorem} \label{theorem-cone}
There exists a rational polyhedral cone $U \subset \Mov{X}^e$ with the
following properties:

(1) $\Pic^0(X_\eta) \cdot U$ intersects the interior of every nef
cone in $\Mov{X}^e$.

(2) $U$ is contained in a union $\cup_{i=1,\ldots,n} \Nef{X_i}^e$ of
  finitely many nef cones in $\Mov{X}^e$.
\end{theorem}
{\bf Proof:} The first step is to choose a rational polyhedral cone
$K_0$ in $\Mov{X}^e$ which maps onto the cone $K$ from the previous
lemma. This is possible by Lemma \ref{lemma-surjective}. More
precisely, for each extremal ray $R_i$ of the cone $K$, choose a
integral vector $w_i$ spanning $R_i$, and a Cartier divisor class $D_i
\in \Mov{X}^e$ such that $p(v_i)=w_i$. This gives a set
$\{D_1,\ldots,D_p\}$ of vectors in $\Mov{X}^e$ spanning a cone which
surjects onto $K$. Define $K_0$ to be the cone spanned by the set
$\{D_1,\cdots,D_p\}$. 

In choosing generators of the cone $K_0$, we are free to replace $D_i$
by $D_i + \Delta$, where $\Delta$ is any element of $T(X/S)$, as long
as the resulting vector still lies in $\Mov{X}^e$. We can use this to
choose generators $D_1,\ldots,D_n$ lying in the big cone. I claim that
if $D \in \Mov{X}^e$ is a class not in the subspace $T(X/S)$, and $A$
is a general very ample divisor on $S$, then $D+f^*(A)$ is big. We use
the numerical criterion for bigness of nef divisors \cite[Proposition
  2.61]{KollarMori1998}: if $N$ is a nef divisor on a $k$-dimensional
projective variety $Y$, then $N$ is big if and only if $N^k
>0$. Choose an SQM $X'$ of $X$ on which $D$ is nef; this implies that
$D+f^*(A)$ is also nef for any ample $A$ on $S$. So consider the
self-intersection number $(D+f^*(A))^3 = D^3+(f^*(A))^3+3D^2 \cdot
f^*(A)+3D \cdot (f^*(A))^2$. The first 2 terms are nonnegative since
$D$ and $f^*(A)$ are nef, and the third term is nonnegative since $D^2
\in \Curv{X'}$ and $f^*(A)$ is nef. I claim the last term must be
strictly positive. If we choose $A$ to be a general very ample divisor
on $S$ , then $A^2$ is a union of $r>0$ points on $S$, and $f^*(A)^2$
is a union of $r$ general fibres of $f$. So $D \cdot (f^*(A))^2 =0$ if
and only if $D \cdot F =0$ for $F$ the class of a fibre of $f$, which
(since $D$ is nef) implies $D \cdot C =0$ for any fibral curve of
$f$. But then $D$ belongs to $T(X'/S)$, which by Lemma
\ref{lemma-numtrivial} equals $T(X/S)$. This contradicts our
assumption on $D$. We conclude that we can add suitable pullbacks of
very ample divisors on $S$ to our generators $D_i$ to ensure they are
big. The benefit of this change is the following: by Theorem
\ref{theorem-localfiniteness}, the decomposition of $\Mov{X}^e$ into
nef cones is now locally finite near the cone $K_0$. In particular
$K_0$ is contained in a union $\cup_{i=1,\ldots,n} \Nef{X_i}^e$ of
finitely many nef cones of SQMs of $X$.

Now choose an SQM $X'$ of $X$ and a Cartier divisor $D$ in the ample
cone $A(X')$. By Lemma \ref{lemma_funddomain} there exists a divisor
$D_0$ in $K_0$ and an element $\phi \in \Pic^0(X_\eta)$ such that
$\phi_*([D_0]) = [D]$ in $N^1(X/S)$.  Therefore in $N^1(X)$ we have
$\phi_*(D_0)=D+\Delta$, where $\Delta$ is a Cartier divisor in
$T(X/S)$. From this expression we cannot conclude that $\phi_*(D_0)$
is ample on $X'$, so the cone $K_0$ is not big enough for our
purposes. 

To remedy this, fix an ample Cartier divisor class $A$ on $S$. For any
$m$ the divisor $f^*(mA)$ is fixed by $\Pic^0(X_\eta)$, so from the
expression above we get $\phi_*(D_0+f^*(mA))=D+\Delta+f^*(mA)$. Now $D$ is
ample on $X'$ and $\Delta$ is numerically trivial on fibral curves,
therefore $D+\Delta$ is ample over $S$. But then by \cite[Proposition
  1.45]{KollarMori1998} $D+\Delta+f^*(mA)=\phi_*(D_0+f^*(mA))$ is ample on
$X'$ for $m$ a sufficiently large positive integer.

Now define $U$ to be the cone spanned by $K_0$ and $f^*(A)$. By the
previous paragraph, $\Pic^0(X_\eta) \cdot U$ intersects the interior
of every nef cone inside $\Mov{X}^e$. It remains to show that $U$ is
contained in a union of finitely many nef cones. To see this, observe
that $f^*(A)$ belongs to every nef cone in $\Mov{X}^e$ since it is
semiample on any SQM of $X$. Therefore $U$ is contained in the same
finite union $\cup_{i=1,\ldots,n} \Nef{X_i}^e$ of nef cones as
  $K_0$. \quad QED

The theorem also proves the full movable cone conjecture for any $X$
satisfying Assumption 1 with the extra assumption that any SQM $X'$ of
$X$ has rational polyhedral nef cone. To see this, note that the extra
assumption implies that the union $\cup_{i=1,\ldots,n} \Nef{X_i}^e$ is
a rational polyhedral cone whose translates by the Mordell--Weil group
cover the entire movable cone. A theorem of Looijenga
\cite[Application 4.15]{Looijenga2009} then says that the existence of
such a cone implies that $\Mov{X}^e$ has a rational polyhedral
fundamental domain for the action of any of the groups
$\Pic^0(X_\eta)$, $\PsAut(X/S)$, or $\PsAut(X,\Delta)$. In practice,
the condition on the nef cones of all SQMs of a given variety may be
difficult to check. However, it has been verified in at least one
nontrivial class of examples, which we now describe.

\begin{examples} \label{example-quadrics}
\textup{Let $X$ be the blowup of $\P^3$ in the 8 basepoints of a
  general net of quadrics. In this case, the anticanonical contraction
  morphism is given by the line bundle $-\frac{1}{2} K_X$ and has base
  $\P^2$. The elliptic fibration $X \arrow \P^2$ has exactly $28=$ $8
  \choose 2$ reducible fibres, corresponding to the set of lines
  through 2 of the 8 basepoints of the net. Each reducible fibre is a
  union $C_1 \cup C_2$ of 2 rational curves meeting transversely in 2
  points, in other words type $I_2$ in Kodaira's classification. (This
  example can therefore be regarded as a global version of Kawamata's
  example \cite[Example 3.8]{Kawamata1997} in which one considers a
  versal deformation of an $I_2$ fibre.)  The SQMs of $X$ are
  sequences of flops in which we first flop a component $C_1$ of a
  reducible fibre, then flop the proper transform $C_2' $ of $C_2$,
  then flop the proper transform $C_1''$ of the other component of the
  fibre containing $C_2'$, and so on. (More generally one can have
  such a sequence of flops for each fibre, but these sequences do not
  interact.) By explicitly computing the action of the Mordell--Weil
  group on certain movable classes, one can prove in this case that
  all the resulting SQMs have rational polyhedral nef cone
  \cite[Theorem 3.5]{Prendergast-Smith2009}.  As explained above, the
  full movable cone conjecture therefore holds for $X$.}
\end{examples}


\section{The nef cone conjecture} \label{section-nefcone}
In this section we collect some results which give evidence for the
nef cone conjecture for our class of varieties. First we observe that
Kawamata's result on the relative nef cone implies that the
$K$-trivial part of the cone of curves obeys the conjecture, and then
we prove finiteness of most types (in the sense of Mori's list) of
divisorial $K$-negative extremal rays. In a slightly different
direction, in the smooth case we use the classification results of
Bauer--Peternell for 3-folds with nef anticanonical bundle to show
that if $X$ is smooth and not rationally connected and $-K_X$ has
Iitaka dimension 2, then the nef cone conjecture again holds for $X$.

As always, we assume throughout this section that $X$ is a
$\Q$-factorial terminal Gorenstein 3-fold with $-K_X$ semiample of
positive Iitaka dimension. The main result (Theorem
\ref{corollary-nefcone}) requires the stronger assumption that $X$ is
smooth and $-K_X$ has Iitaka dimension 2.

\begin{proposition}
$\Curv{X} \cap K^\perp = \Curv{X/S}$.
\end{proposition}
{\bf Proof:} The inclusion $\Curv{X/S} \subset \Curv{X} \cap K^\perp$
is immediate since $f:X \arrow S$ is given by a power of
$-K_X$ ; it remains to prove the reverse inclusion. Suppose it does
not hold: then we can find a point $x \in \Curv{X} \cap K^\perp$ which
does not belong to $\Curv{X/S}$. Choose a rational class $H \in
N^1(X)$ such that $H$ is strictly positive on $\Curv{X/S} \backslash
\{0\}$, but $H \cdot x < 0$. By the Kleiman condition, $H$ is ample
over $S$. By \cite[Proposition 1.45]{KollarMori1998}, if $A$ is any
ample class on $S$, then $H+f^*(nA)$ is ample on $X$ for $n$ a
sufficiently large natural number.

Now the morphism $f: X \arrow S$ is given by $-mK_X$ for some positive
integer $m$, so by the contraction theorem there is an ample line
bundle $A$ on $S$ such that $f^*(A) = -mK_X$. Putting this together
with the previous paragraph, we see that $H-nK_X$ is ample on $X$ for
$n$ a sufficiently large and divisible positive integer. But since by
assumption $x$ is a point in $K_X^\perp$, we have $(H-nK_X) \cdot x =
H \cdot x <0$, which since $H-nK_X$ is ample contradicts the
assumption $x \in \Curv{X}$. \quad QED


\begin{corollary} \label{corollary-finiteness}
$\Curv{X}$ has only finitely many extremal rays in $K^\perp$
  corresponding to divisorial contractions or fibre space structures,
  up to the action of $\Aut(X/S)$.
\end{corollary} {\bf Proof:} By the previous lemma, an extremal ray of
$\Curv{X}$ lying in $K^\perp$ must be an extremal ray of $\Curv{X/S}$;
moreover, if such an extremal ray corresponds to a contraction
morphism, then that contraction is over $S$. Now Kawamata's theorem on
the relative nef cone \cite[Theorem 3.6, Theorem 4.4]{Kawamata1997}
shows that the number of contractions over $S$ is finite up to the
action of $\Aut(X/S)$.\quad QED

Next we consider the $K$-negative extremal rays of $\Curv{X}$. By the
cone theorem any such ray can be contracted, and so corresponds either
to a Mori fibre space structure of $X$ or to a divisorial contraction
of one of the 5 types on Mori's list (Theorem
\ref{theorem-moriclassification}). The next result shows that with the
exception of one type, the divisorial contractions of $X$ must always
be finite in number.

\begin{proposition}
The number of $K$-negative extremal contractions of $X$ of types 2, 3,
4 and 5 on Mori's list is finite.
\end{proposition}
{\bf Proof:} I claim that if $\phi: X \arrow Y$ is a $K$-negative
extremal contraction of type 2, 3, 4 or 5 on Mori's list, the
exceptional divisor $D_\phi$ of $\phi$ must be disjoint from the
exceptional divisor of any other $K$-negative extremal
contraction. Given this, the proposition follows, since otherwise we
could find a morphism $f: X \arrow Z$ of projective varieties
contracting an arbitrary number (in particular more than $\rho(X)-1$)
of prime divisors on $X$ whose classes in $N^1(X)$ are linearly
independent.

To prove the claim, suppose $\phi: X \arrow Y$ is of type 2, 3, 4 or
5, and $\psi: X \arrow Z$ is any $K$-negative extremal contraction. If
$D_\phi$ is not disjoint from the exceptional divisor $D_\psi$ of
$\psi$, then they intersect in a 1-dimensional cycle supported on a
curve $C$. So $\psi$ contracts the curve $C \subset D_\phi$. But the
description in Theorem \ref{theorem-moriclassification} of the
exceptional divisor $D_\phi$ tells us that for $\phi$ of type 2, 3, 4
or 5, all curves in $D_\phi$ are numerically proportional in $N_1(X)$,
so $\psi$ must contract all of $D_\phi$. Since the exceptional divisor
of $\psi$ is prime, we conclude that $\psi = \phi$. \quad QED

To summarise, we have shown finiteness up to automorphisms of the
number of $K$-trivial extremal contractions, and finiteness of the
number of divisorial $K$-negative extremal contractions of all but one
type on Mori's list. To complete the proof of the nef cone conjecture
for the given class of varieties, it remains to show finiteness up to
automorphisms of the set of divisorial contractions of type 1 on
Mori's list and the set of $K$-negative fibre space structures of
$X$. I am unable to do this at present.


We conclude by mentioning the classification results of
Bauer--Peternell for smooth 3-folds with $-K_X$ nef, and
showing how these reduce the nef cone conjecture to the rationally
connected case. We should point out that the theorem below is only a
small part of Bauer--Peternell's results, which give information about
all smooth 3-folds with $-K_X$ nef and numerically nontrivial. See
\cite{BauerPeternell2004} for details. From now on, we consider only
smooth $X$ with $-K_X$ semiample of Iitaka dimension 2.

\begin{theorem}[Bauer--Peternell] \label{theorem-bpclass}
Suppose that $X$ is a smooth 3-fold with $-K_X$ semiample of
Iitaka dimension 2. Then one of the following holds:
\begin{enumerate}
\item $X$ is rationally connected;
\item $X \iso_f S \times B$ with $S$ a del Pezzo surface, $B$ an
  elliptic curve;
\item $X \iso_f \P(E) \times B$ with $E$ a rank-2 bundle over an
  elliptic curve $A$, and $B$ another elliptic curve.
\end{enumerate}
Here the symbol $\iso_f$ means that some finite etale cover of $X$ is
isomorphic to the given variety. 
\end{theorem}
Bauer--Peternell give analogous results in the case $-K_X$ of Iitaka
dimension 1, but I was not able to apply these results to the cone conjecture.

\begin{corollary} \label{corollary-bp}
For $X$ as in Theorem \ref{theorem-bpclass}, if $X$ is not rationally
connected, then $\Curv{X}$ has finitely many $K$-negative extremal
rays.
\end{corollary} {\bf Proof of Corollary:} By Theorem \ref{theorem-bpclass} if
$X$ is not rationally connected there is a finite etale cover $\pi: Y
\arrow X$ such that $Y$ has the one of the two product structures
listed. I claim that the set of $K$-negative extremal rays of
$\Curv{X}$ can be identified naturally with a subset of the
$K$-negative extremal rays of $\Curv{Y}$, and that the latter set is
finite.

First note that since $\pi$ is etale, the sheaf of relative
differentials $\Omega_{Y/X}$ is zero, so the exact sequence $\pi^*
\Omega_X \arrow \Omega_Y \arrow \Omega_{Y/X} \arrow 0$ implies that
$\Omega_Y = \pi^* \Omega_X$, which in turn gives $K_Y = \pi^* K_X$. In
particular, if $C$ is any curve on $X$, we have $K_X \cdot C <0$ if
and only if $K_Y \cdot \pi^{-1}(C) <0$. Now suppose $R$ is a
$K$-negative extremal ray of $\Curv{X}$. By the cone theorem, there is
a semiample line bundle $L$ on $X$ such that $\Curv{X} \cap L^\perp =
R$. The pullback $\pi^* L$ is again semiample, and clearly contracts
no curves except those whose image on $X$ have class in the ray
$R$. This proves the first claim, that the set of $K$-negative
extremal rays of $\Curv{X}$ can be identified with a subset of the
$K$-negative extremal rays of $\Curv{Y}$. It remains to show that
$\Curv{Y}$ has only finitely many $K$-negative extremal rays. This
follows easily from the product structure of $Y$, as we now explain.

Case 1: here $Y = S \times B$, with $S$ a del Pezzo surface, $B$ an
elliptic curve. Denote the projections of $Y$ onto its two factors by
$pr_1$ and $pr_2$. Suppose $R$ is a $K$-negative extremal ray of
$\Curv{Y}$. Then $R$ is spanned by the class of a rational curve
$C$. There cannot be a nonconstant morphism $C \arrow B$ since $B$ has
genus 1, so $C$ must be contained in a surface $pr_2^{-1}(p) \iso S$
inside $Y$, where $p \in B$ is a point. Therefore the ray $R$ must
belong to the subcone $\Curv{Y/B} \subset \Curv{Y}$. Since $Y$ is a
product $S \times B$, we have $\Curv{Y/B} \iso \Curv{S}$, which is
rational polyhedral since $S$ is del Pezzo. So there are only finitely
many possibilities for $R$.

Case 2: here $Y = \P(E) \times B$, with $E$ a rank-2 bundle on an
elliptic curve $A$, and $B$ another elliptic curve. As above, any
rational curve on $Y$ must project to a point on $B$, so any
$K$-negative extremal ray $R \subset \Curv{Y}$ must belong to the
subcone $\Curv{Y/B}$. Again this is just the cone $\Curv{\P(E)}$. Any
ruled surface has Picard number 2, so its closed cone of curves has
exactly 2 extremal rays. \quad QED

By Proposition \ref{corollary-finiteness} the set of extremal rays in
$K^\perp$ corresponding to birational contractions or fibre space
structures is in any case finite up to the action of $\Aut(X/S)$, so we
can put these results together to get the following.

\begin{theorem} \label{corollary-nefcone}
Suppose $X$ is a smooth 3-fold with $-K_X$ semiample of Iitaka
dimension 2, and $X$ is not rationally connected. Then $\Curv{X}$ has
finitely many extremal rays corresponding to birational contractions
or fibre space structures, up to the action of $\Aut(X)$.
\end{theorem}


\small
\sc

Leibniz Universit\"at Hannover, Institut f\"ur Algebraische Geometrie,
Welfengarten 1, D-30167 Hannover, Germany

{\it Email address:} {\tt artie@math.uni-hannover.de}
\end{document}